\documentstyle[bezier,12pt]{article}

\input{emlines.sty}

\newenvironment{slshape}{\sl}{}
\newcommand{\mathcal}[1]{{\cal #1}}
\newcommand{\matrm}[1]{\mbox{#1}}
\newcommand{\emph}[1]{{\em #1\/}}

\newtheorem{thm}{Theorem}[section]
\newtheorem{prop}[thm]{Proposition}
\newtheorem{lem}[thm]{Lemma}
\newtheorem{corol}[thm]{Corollary}
\newenvironment{theorem}{\begin{thm}\begin{slshape}}{\end{slshape}\end{thm}}
\newenvironment{proposition}{\begin{prop}\begin{slshape}}{\end{slshape}\end{prop}}
\newenvironment{lemma}{\begin{lem}\begin{slshape}}{\end{slshape}\end{lem}}
\newenvironment{corollary}{\begin{corol}\begin{slshape}}{\end{slshape}\end{corol}}
\newtheorem{dfn}[thm]{Definition}
\newenvironment{definition}{\begin{dfn}\em}{\end{dfn}}
\newtheorem{rem}[thm]{Remark}
\newenvironment{remark}{\begin{rem}\em}{\end{rem}}
\newtheorem{que}[thm]{Question}
\newenvironment{question}{\begin{que}\em}{\end{que}}
\newtheorem{xmp}[thm]{Example}

\newenvironment{claim}[1]{{\par\bigbreak\noindent\sl Claim #1. }}{\par\bigbreak}
\newcommand{\phase}[1]{\par\smallbreak\noindent{\sc Phase #1}}
\newcommand{\condition}[1]{\par\hangindent=\parindent
\noindent{\sl Condition #1. }}
\newcommand{\itemm}[1]{\par\smallbreak\noindent{\sl Item #1. }}
\newcommand{\convention}[1]{\par\medbreak\noindent{\sl Convention. }#1
\par\medbreak}
\newcommand{\notation}[1]{\par\medbreak\noindent{\sl Notation. }#1
\par\medbreak}
\newcommand{\comment}[2]{\par\medbreak\noindent {\sl Comment. }#2
\par\medbreak}
\newcommand{\case}[2]{\par\bigbreak\noindent{\sl Case #1}: #2.\par\noindent}
\newcommand{\subcase}[2]{\par\medbreak\noindent{\sl Subcase #1}: #2.\par\noindent}
\newcommand{\rrule}[1]{\par\hangindent=\parindent
\noindent{\sl Rule #1. }}

\newcommand{\qed}{$\;\;\;\Box$}
\newcommand{\bull}{$\;\;\;$\vrule height .9ex width .8ex depth -.1ex}
\newenvironment{proof}{\par\smallbreak\noindent{\sl Proof.~}}
{\unskip\nobreak\hfill \bull \par\medbreak}

\newenvironment{subproof}{\par\smallbreak\noindent{\sl Proof.~}}
{\unskip\nobreak\hfill \qed \par\medbreak}

\newcommand{\one}{$\cal{A}$}
\newcommand{\two}{$\cal{B}$}
\newcommand{\oneo}{${\cal A}_0$}
\newcommand{\twoo}{${\cal B}_0$}
\newcommand{\oneef}{spoiler}
\newcommand{\twoef}{duplicator}

\newcommand{\V}{{\cal V}}
\newcommand{\E}{{\cal E}}
\newcommand{\lineg}[1]{{\cal L}( #1 )}
\newcommand{\EF}{Ehrenfeucht-Fra\"\i ss\'e}
\newsavebox{\symbox}
\savebox{\symbox}{sym}
\newsavebox{\symplusbox}
\savebox{\symplusbox}{sym\raisebox{0.1ex}{$\scriptstyle{}+{}$}}
\newsavebox{\efbox}
\savebox{\efbox}{\scriptsize EF}
\newcommand{\lsym}[1]{L_{\usebox{\symbox}} ( #1 )}
\newcommand{\lsymplus}[1]{L_{\usebox{\symplusbox}} ( #1 )}
\newcommand{\llsym}[1]{L'_{\usebox{\symbox}} ( #1 )}
\newcommand{\lef}[1]{L_{\usebox{\efbox}} ( #1 )}
\newcommand{\avoid}[2]{\mbox{AVOID}( #1, #2)}
\newcommand{\gamesym}[1]{\mbox{SYM}( #1 )}
\newcommand{\symr}[2]{\mbox{SYM}_{#1}( #2 )}
\newcommand{\symplus}[1]{\mbox{SYM\raisebox{0.3ex}{$+$}}( #1 )}
\newcommand{\efgame}[1]{\mbox{EF}( #1 )}
\newsavebox{\bsymm}
\savebox{\bsymm}{$\mathcal{C}_{\matrm{sym}}$}
\newcommand{\symm}{\usebox{\bsymm}}

\newcommand{\refeq}[1]{(\ref{#1})}
\newcommand{\function}[2]{:#1 \rightarrow #2}
\newcommand{\of}[1]{\left( #1 \right)}
\newcommand{\eqdef}{\stackrel{\rm def}{=}}

\newcommand{\And}{\wedge}
\newcommand{\Or}{\vee}

\title{On the Lengths of Symmetry\\
Breaking-Preserving Games on Graphs}

\author{Frank~Harary\thanks{Computer Science Department,
New Mexico State University,
Las Cruces, NM 88003, USA.}\quad
Wolfgang~Slany\thanks{Institut f\"ur Informationssysteme, Technische Universit\"at Wien,
Favoritenstr.~9, A-1040 Wien, Austria. Research partly supported from
Austrian Science Foundation grant Z29-INF.}\quad
Oleg~Verbitsky\thanks{
Department of Mechanics \& Mathematics, Lviv University,
Universytetska~1, 79000 Lviv, Ukraine.}}

\date{}

\begin{document}
\maketitle

\begin{abstract}
Given a graph $G$, we consider a game where two players, \one\/
and \two, alternatingly color edges
of $G$ in red and in blue respectively.
Let $\lsym G$ be the maximum number of moves in which \two\/
is able to keep the red and the blue subgraphs isomorphic,
if \one\/ plays optimally to destroy the isomorphism.
This value is a lower bound for the duration of any {\em avoidance\/}
game on $G$ under the assumption that \two\/ plays optimally.
We prove that if $G$ is a path or a cycle of odd length $n$, then
$\Omega(\log n)\le \lsym G\le O(\log^2 n)$. The lower bound
is based on relations with \EF\/ games from model theory.
We also consider complete graphs and prove that $\lsym{K_n}=O(1)$.
\end{abstract}

\section{Introduction}

The symmetry breaking-preserving game $\gamesym G$ is played by
two players on a graph $G$. The players, \one\/ and \two, alternatingly
color edges of $G$ in red and in blue respectively, one edge per move.
Player \one\/ is first to move. A round of the game consists of
a move of \one\/ and the following move of \two. The objective of
\two\/ is to keep the red and the blue subgraphs of $G$ isomorphic
after every round. As soon as \two\/ fails to do so, this is a win
for \one. If \two\/ succeeds until all the edges are colored,
this is a win for him.

This game was introduced in \cite{HSV} in the context of the
graph avoidance games \cite{Har,EHa}. The game $\avoid GF$ is a two-person
edge-coloring game on a graph $G$ with the following ending
condition: The player who first creates a monochromatic copy of
a forbidden subgraph $F$ loses. As easily seen, as long as \two\/
does not lose in $\gamesym G$, he does not lose in $\avoid GF$
for any $F$.

In \cite{HSV} we addressed the class $\symm$ of those graphs $G$
for which \two\/ has a winning strategy in $\gamesym G$.
We now consider a more general problem:
Given $G$, how long is \two\/ able to keep
the red and the blue subgraphs isomorphic if both players play optimally?
We define $\lsym G$, the \emph{length} of the game $\gamesym G$,
to be the maximum number of rounds in which \two, playing optimally,
does not lose, independently of \one's strategy
(a precise definition is given in Section \ref{s:defn}).
This function of a graph $G$ will be our main concern.

Note that $G$ belongs to $\symm$ iff $\lsym G$ has the maximum
possible value $\lfloor m/2\rfloor$,
where $m$ is the size of $G$. In \cite{HSV} we observe that
$\symm$ contains all graphs having an involutory automorphism
without fixed edges. Though generally it is NP hard to recognize
if such an automorphism exists for a given graph, we nevertheless
obtain many examples of graphs in $\symm$. The simplest examples
are even cycles and paths. Let $C_n$ (resp.\ $P_n$) denote the cycle
(resp.\ path) of size $n$. Thus, we have
$\lsym{C_n}=\lsym{P_n}=n/2$ for $n$ even.

In the present paper we treat odd cycles and paths.
For odd $n$, we prove that
$$
\Omega(\log n)\le \lsym{P_n}\le O(\log^2 n), \quad
\Omega(\log n)\le \lsym{C_n}\le O(\log^2 n).
$$
Our proof of the lower bound is based on the connections with
\EF\/ games known in model theory \cite{fmt}, what may be of independent
interest. In particular, we use the well known fact that
the length of the \EF\/ game on the pair of paths $P_n$ and $P_{n+1}$
equals $\log n$ up to an additive constant and the same is true for
the pair of cycles $C_n$ and $C_{n+1}$.

We also consider symmetry breaking-preserving games on complete graphs.
As implicitly shown in \cite{HSV}, $\lsym{K_n}\le n-2$.
We now improve this estimate showing that
$\lsym{K_n}$ is for all $n$ bounded by an absolute constant.

Note that all the upper (resp.\ lower) bounds proven here are
based on efficiently computable strategies for the player \one\/
(resp.\ \two).

In the next section we give the definitions and state some useful facts.
We estimate the asymptotics of $\lsym G$ for odd paths and
cycles in Section \ref{s:odd} and for complete graphs in Section \ref{s:k_n}.

\section{Preliminaries}\label{s:defn}

Given a graph $G$, we denote its vertex set by $\V(G)$ and its edge set
by $\E(G)$.

The \emph{symmetry breaking-preserving game} on a graph $G$,
denoted by $\gamesym G$, is a two-person positional game of
the following kind.
Two players, \one\/ and \two, alternatingly color edges of a graph
$G$ in red and in blue respectively. Player \one\/ starts the game.
In a {\em move}, a player colors an edge that was so far uncolored.
The $i$-th {\em round\/} consists of the $i$-th move of \one\/
and the $i$-th move of \two.
Let $a_i$ (resp.\ $b_i$) denote an edge
colored by \one\/ (resp.\ \two) in the $i$-th round.
Let $A_i=\{a_1,\ldots,a_i\}$ (resp.\ $B_i=\{b_1,\ldots,b_i\}$) consist of
the red (resp.\ blue) edges colored up to the $i$-th round.
Player \two\/ wins in $\gamesym G$ if the subgraphs $A_i$ and $B_i$ are
isomorphic for every $i\le |\E(G)|/2$. As soon as an isomorphism
between $A_i$ and $B_i$ is violated, this is a win for \one.

A {\em strategy\/} for a player determines the edge to be colored by him at every
round of the game. Formally, let $\epsilon$ denote the empty sequence.
A strategy of \one\/ is a function $S_1$ that maps every, possibly empty,
sequence of pairwise distinct edges $e_1,\ldots,e_i$
into an edge different from $e_1$, \ldots, and $e_i$ and from
$S_1(\epsilon)$, $S_1(e_1)$,
$S_1(e_1,e_2)$, \ldots, and $S_1(e_1,\ldots,e_{i-1})$.
A strategy of \two\/ is a function $S_2$ that maps every nonempty
sequence of pairwise distinct edges $e_1,\ldots,e_i$
into an edge different from $e_1$, \ldots, and $e_i$ and from
$S_2(e_1)$, $S_2(e_1,e_2)$, \ldots, and $S_2(e_1,\ldots,e_{i-1})$.
If \one\/ follows a strategy $S_1$ and \two\/ follows a strategy $S_2$,
then $a_i=S_1(b_1,\ldots,b_{i-1})$ and $b_i=S_2(a_1,\ldots,a_{i})$.

The {\em length\/} of the game is the total number of rounds under the
condition that the players play optimally. To be more precise,
assume that \one\/ follows a strategy $S_1$ and \two\/ follows
a strategy $S_2$ and let $l(S_1,S_2)$ denote the maximum $l$ such that
$A_i$ and $B_i$ are isomorphic for every $i\le l$.
We denote the length of $\gamesym G$ by $\lsym G$ and define it by
$$
\lsym G=\max_{S_2}\min_{S_1} l(S_1,S_2).
$$
An alternative definition could be
$$
\llsym G=\min_{S_1}\max_{S_2} l(S_1,S_2).
$$
Observe that the definitions are equivalent.
\begin{proposition}\label{prop:llsym}
$\lsym G = \llsym G$.
\end{proposition}
\begin{proof}
The inequality $\lsym G\le \llsym G$ holds true by the universal
min-max relation. To prove the reverse inequality, define a game
$\symr r G$ to be a variant of $\gamesym G$ in which \two\/ wins if
he does not lose the first $r$ rounds. A strategy of a player is
\emph{winning} if it beats every strategy of the opponent.
Since $\symr r G$ is a finite perfect information game with no draws,
in this game one of the players has a winning strategy.
Assume that $\lsym G=l$ and $l<\lfloor|\E(G)|/2\rfloor$.
This means that \two\/ has no winning strategy in $\symr{l+1}G$.
Hence the winning strategy in $\symr{l+1}G$ exists for \one,
which implies that $\llsym G\le l$.
\end{proof}

We will refer to the following observation that follows from~\cite{HSV}.

\begin{proposition}\label{prop:auto}
If $G$ has an involutory automorphism without fixed edges, then
$$
\lsym G = |\E(G)|/2.
$$
\end{proposition}

\begin{proof}
If $\phi\function{\V(G)}{\V(G)}$ is an automorphism of $G$, it
determines a permutation $\phi'\function{\E(G)}{\E(G)}$ by
$\phi'(\{u,v\})=\{\phi(u),\phi(v)\}$. We assume that $\phi$
is involutory and $\phi'$ has no fixed element. Then
the edge set $\E(G)$ is partitioned into 2-subsets of the form
$\{e,\phi'(e)\}$. This gives \two\/ the following winning strategy:
Whenever \one\/ chooses an edge $e$, \two\/ chooses
the edge $\phi'(e)$. After every round of the game, an isomorphism
between the red and the blue subgraphs is induced by $\phi$.
\end{proof}

The \emph{\EF\/ game} can be played on an arbitrary structure.
We give a definition conformably to graphs.
Assume that graphs $G_0$ and $G_1$ have disjoint vertex sets.
In the Ehren\-feucht-Fra\"\i ss\'e game on $G_0$ and $G_1$,
denoted further on by $\efgame{G_0,G_1}$, the players \one\/
and \two\/ alternatingly pick up vertices of either $G_0$ or $G_1$,
one vertex per move. \one\/ starts the game. Let $u_i$ (resp.\ $v_i$)
be the vertex picked up by \one\/ (resp.\ by \two) in his $i$-th
move. In each round the objective of \two\/ is to obey the
following conditions.
\begin{itemize}
\item
If $u_i\in\V(G_a)$, then $v_i\in\V(G_{1-a})$.
\item
The correspondence ``$u_i$ to $v_i$'' is a partial isomorphism
between $G_0$ and $G_1$, i.e., an isomorphism between the subgraphs
of $G_0$ and $G_1$ induced by the chosen vertices.
\end{itemize}
The maximum number of rounds in which \two\/,
irrespective of \one's strategy, is able to obey these two conditions
is denoted by $\lef{G_0,G_1}$ and is formally defined
similarly to $\lsym G$.

Let $\log n$ denote logarithm base 2.
We will use the following folklore result.

\begin{proposition}\label{prop:lef}
For every $n$,
\begin{enumerate}
\item
$\log n-2<\lef{P_n,P_{n+1}}<\log n+2$.
\item
$\log n-1<\lef{C_n,C_{n+1}}<\log n+1$.
\end{enumerate}
\end{proposition}

The proof can be found in \cite[Theorems 2.1.2 and 2.1.3]{Spe}
for the case of paths. The case of cycles can be treated similarly
(cf.\ \cite[Example 2.3.8]{fmt}).

\section{Games on paths and cycles}\label{s:odd}

Given two functions
$f(n)$ and $g(n)$, we use notation $f(n)=\Omega(g(n))$ whenever
$f(n)\ge c\cdot g(n)$ for some $c>0$ and all $n$.

The main result of this section estimates the asymptotics of
$\lsym G$ for odd paths and cycles. It should be contrasted with
even paths and cycles, for which $\lsym{P_n}=\lsym{C_n}=n/2$
by Proposition \ref{prop:auto}.

\begin{theorem}\label{thm:odd}
If $n$ is odd, then
\begin{enumerate}
\item
$\lsym{P_n}=\Omega(\log n)$ and
$\lsym{C_n}=\Omega(\log n)$,
\item
$\lsym{P_n}=O(\log^2 n)$ and
$\lsym{C_n}=O(\log^2 n)$.
\end{enumerate}
\end{theorem}
The proof of the theorem is given in the rest of this section.

\subsection{Lower bound}\label{ss:lower}

To prove the lower bounds for $\lsym{P_n}$ and $\lsym{C_n}$, we
relate the symmetry breaking-preserving game with the \EF\/ game.
We are actually able to prove the claim 1 of Theorem \ref{thm:odd}
in two different ways, both using Proposition \ref{prop:lef}.
We will refer to one way as the logical approach and to the other
way as the combinatorial approach.

We start with the brief overview of the logical approach for the
case of paths; the case of cycles is virtually identical.
Given an odd path $P_n$, we consider also the even path $P_{n+1}$
for which we know that $\lsym{P_{n+1}}=(n+1)/2$. As known from
model theory, the length of $\efgame{P_n,P_{n+1}}$ tells us to which
extent the properties of $P_{n+1}$ expressible in first order logic
hold true for $P_n$ (see Lemma \ref{lem:ef}). On the other hand,
Lemma \ref{lem:loglsym} tells us to which extent the property of
a graph $G$ that $\lsym G\ge k$ is first order expressible.
Putting it together, we see that, as the property $\lsym G=\Omega(\log n)$
is true for $P_{n+1}$, it must be true for $P_n$ too
(cf.\ Proposition \ref{prop:sym_ef}).
Curiously, this method proves the existence of the desired strategy for
player \two\/ without yielding it explicitly.

The combinatorial approach does not exploit the logical aspects
of the \EF\/ game. Instead, it directly exploits the partial isomorphism
constructed during the course of $\efgame{P_n,P_{n-1}}$ in order to
translate, as long as possible, the winning strategy of \two\/ from
$\gamesym{P_{n+1}}$ into $\gamesym{P_n}$ (see Proposition~\ref{prop:new}).

The combinatorial approach gives us a bound twice as good as the
logical approach. What is more important than the gain in a
multiplicative constant, the former approach provides us with
an efficiently computable strategy for \two.
Nevertheless, though the combinatorial approach is more preferable
to the logical one in the particular cases of odd paths and cycles,
generally it has a more restrictive applicability range.
We do not exclude that both techniques may be useful in the
analysis of other games on graphs (cf.\ Remark~\ref{rem:ext}).

We now present both the proof methods in detail, starting from the
logical one.

{}From the logical point of view a graph $G$ is a structure consisting
of a single binary predicate $E$ on $\V(G)$ such that $E(u,v)$ iff
$u$ and $v$ are adjacent. Every closed first order formula over
vocabulary $\{E,{=}\}$ is either true or false on $G$.

\begin{lemma}{\bf (\cite[Theorem 2.2.8]{fmt})}\label{lem:ef}
$G_0$ and $G_1$ satisfy precisely the same first order sentences
with at most $\lef{G_0,G_1}$ quantifiers.
\end{lemma}

\noindent Observe that the sentence ``$\lsym G \ge k$'' is expressible with
$4k$ quantifiers.

\begin{lemma}\label{lem:loglsym}
There is a first order formula $\Phi_k$ with $4k$ quantifiers
that is true on $G$ of size at least $2k$ iff $\lsym G\ge k$.
\end{lemma}

\begin{proof}
Let $DIST(x_1,x_2,y_1,y_2)$ express the property that two pairs of vertices
$\{x_1,x_2\}$ and $\{y_1,y_2\}$ are distinct. Formally,
$$
DIST(x_1,x_2,y_1,y_2)\eqdef
\neg\of{(x_1=y_1\And x_2=y_2)\Or(x_1=y_2\And x_2=y_1)}
$$
Let $u_{1,1},u_{1,2},\ldots,u_{k,1},u_{k,2}$ and
$v_{1,1},v_{1,2},\ldots,v_{k,1},v_{k,2}$ be variables ranging over
$\V(G)$ with meaning that in the $i$-th round \one\/ chooses
an edge $\{u_{i,1},u_{i,2}\}$ and \two\/ chooses an edge
$\{v_{i,1},v_{i,2}\}$. We also need a formula $ISO_j$ to express the
fact that the subgraphs consisting of the edges chosen by the players
during the first $j$ rounds are isomorphic:
\begin{eqnarray*}
ISO_j(u_{1,1},u_{1,2},\ldots,u_{j,1},u_{j,2},
v_{1,1},v_{1,2},\ldots,v_{j,1},v_{j,2})\eqdef \\
\bigvee_{
\begin{array}{c}
\scriptstyle f
\end{array}
}
\bigwedge_{
\begin{array}{c}
\scriptstyle 1\le i,i'\le j \\
\scriptstyle 1\le a,a'\le 2
\end{array}
}
\of{u_{i,a}=u_{i',a'}
\leftrightarrow
v_{f(i,a)}=v_{f(i',a')}
},
\end{eqnarray*}
where the disjunction is over all permutations of the index set
$\{1,\ldots,j\}\times\{1,2\}$ with the property that if
$f(i,1)=(m,a)$, then $f(i,2)=(m,3-a)$ for all $i\le j$.
The permutation $f$ should be thought of as a map from the
multiset $\{u_{i,a}\}_{i\le j; a\le 2}$ to the
multiset $\{v_{i,a}\}_{i\le j; a\le 2}$ taking every edge
$\{u_{i,1},u_{i,2}\}$ to some edge $\{v_{m,1},v_{m,2}\}$.
Such a permutation is a subgraph isomorphism if it takes equal $u$'s
to equal $v$'s and distinct $u$'s to distinct $v$'s.

Define formulas
$$
A_j\eqdef E(u_{j,1},u_{j,2}) \And
\bigwedge_{i=1}^{j-1} DIST(u_{j,1},u_{j,2},u_{i,1},u_{i,2})\And
\bigwedge_{i=1}^{j-1} DIST(u_{j,1},u_{j,2},v_{i,1},v_{i,2})
$$
and
$$
B_j\eqdef E(v_{j,1},v_{j,2})\And
\bigwedge_{i=1}^{j} DIST(v_{j,1},v_{j,2},u_{i,1},u_{i,2})\And
\bigwedge_{i=1}^{j-1} DIST(v_{j,1},v_{j,2},v_{i,1},v_{i,2})
$$
saying that $\{u_{j,1},u_{j,2}\}$ and, respectively, $\{v_{j,1},v_{j,2}\}$
are edges different from the edges chosen by the players previously.
The formula
\begin{eqnarray*}
\Phi_k&\eqdef&
\forall u_{1,1}\forall u_{1,2} \exists v_{1,1}\exists v_{1,2}\ldots
\forall u_{k,1}\forall u_{k,2} \exists v_{k,1}\exists v_{k,2}
\,\Biggl(\,\,\bigwedge_{j=1}^k A_j\longrightarrow\\ &&\quad
\bigwedge_{j=1}^k B_j\And
\bigwedge_{j=1}^k ISO_j(u_{1,1},u_{1,2},\ldots,u_{j,1},u_{j,2},
v_{1,1},v_{1,2},\ldots,v_{j,1},v_{j,2})\Biggr)
\end{eqnarray*}
is as desired.
Indeed, assume that \two\/ has a strategy non-losing $k$ rounds.
Then $\Phi_k$ is true because, if all $u_{j,1}$, $u_{j,2}$ are
chosen so that the antecedent in $\Phi_k$ is satisfied, then
$v_{j,1}$, $v_{j,2}$ satisfying the consequent are provided
by \two's strategy.

On the other hand, if $\Phi_k$ is true, then the following strategy
of \two\/ does not lose $k$ rounds to any strategy of \one.
We describe the $j$-th move of \two. Assume that \one\/ and \two\/
have previously chosen edges $\{u_{1,1}, u_{1,2}\}$, \ldots,
$\{u_{j,1}, u_{j,2}\}$ and $\{v_{1,1}, v_{1,2}\}$, \ldots,
$\{v_{j-1,1}, v_{j-1,2}\}$. In particular, $u_{1,1}, u_{1,2}, \ldots,
u_{j,1}, u_{j,2}$ satisfy the antecedent in $\Phi_k$.
Then \two\/ chooses an edge $\{v_{j,1}, v_{j,2}\}$ with vertices
$v_{j,1}$ and $v_{j,2}$ whose existence is claimed by $\Phi_k$.
Such $v_{j,1}$ and $v_{j,2}$ satisfy the members $B_j$ and $ISO_j$
of the consequent in $\Phi_k$
because otherwise one could choose the subsequent
$u_{j+1,1},u_{j+1,2},\ldots,u_{k,1},u_{k,2}$ satisfying the antecedent
in $\Phi_k$ and therewith falsify the implication. It follows that
this move of \two\/ is legitimate and successful.
\end{proof}

\begin{proposition}{\bf (logical approach)}\label{prop:sym_ef}
$\lsym{G_1}\ge\min\{\frac14\lef{G_0,G_1},\lsym{G_0}\}$
\end{proposition}

\begin{proof}
Assume that $\lsym{G_0}\ge k$ and $\lef{G_0,G_1}\ge 4k$.
The former inequality implies that $G_0$ has size at least $2k$.
By the latter inequality, the same must be also true for $G_1$.
By Lemma \ref{lem:loglsym}, $G_0$ satisfies $\Phi_k$.
By Lemma \ref{lem:ef}, $G_1$ also satisfies $\Phi_k$ and
therefore, again by Lemma \ref{lem:loglsym}, $\lsym{G_1}\ge k$.
\end{proof}

We now turn to the combinatorial approach to the proof
of Theorem~\ref{thm:odd}~(1).

Given a graph $H$, let $\lineg H$ denote its {\em line graph}.
Recall that $\V(\lineg H)=\E(H)$ and two vertices $e_1$ and $e_2$
of $\lineg H$ are connected by an edge in this graph iff they
have a common vertex in $H$. Two graphs $H_1$ and $H_2$ are
{\em edge-isomorphic\/} if there is a one-to-one map from
$\E(H_1)$ onto $\E(H_2)$ preserving the adjacency of edges.
In other words, $H_1$ and $H_2$ are edge-isomorphic iff $\lineg{H_1}$
and $\lineg{H_2}$ are isomorphic. If two graphs are isomorphic,
then they are obviously edge-isomorphic.
The Whitney theorem \cite[Theorem~8.3]{Har1} says that the converse implication
is also true for all connected $H_1$ and $H_2$ unless one of them is $K_3$
and the other $K_{1,3}$.

To avoid ambiguity, in the next proposition we keep the names
\one\/ and \two\/ for the players in the game $\gamesym{G_1}$,
but rename them \oneo\/ and \twoo\/ in $\gamesym{G_0}$,
and {\em \oneef\/} and {\em \twoef\/} in $\efgame{G_0,G_1}$.

\begin{proposition}{\bf (combinatorial approach)}\label{prop:new}
If $G_1$ does not contain a subgraph $K_3$, then
$\lsym{G_1}\ge\min\{\frac12\lef{\lineg{G_0},\lineg{G_1}},\lsym{G_0}\}$.
Moreover, the player \two\/ in $\gamesym{G_1}$ has an efficiently
computable strategy $S$ with oracle access to a strategy $D$ of
the \twoef\/ in $\efgame{\lineg{G_0},\lineg{G_1}}$ and to a strategy $S_0$
of \twoo\/ in $\gamesym{G_0}$ such that, if $D$ does not lose $l$
rounds irrespective of the \oneef's strategy and $S_0$ does not lose
$m$ rounds irrespective of \oneo's strategy, then $S(D,S_0)$ does not lose
at least $\min\{\frac12l,m\}$ rounds irrespective of \one's strategy.
\end{proposition}

\begin{proof}
To make a move according to $S(D,S_0)$, in each round of $\gamesym{G_1}$
the player \two\/ simulates one round of $\gamesym{G_0}$ following $S_0$
and two rounds of $\efgame{\lineg{G_0},\lineg{G_1}}$ following $D$.
Before describing $S(D,S_0)$, we introduce some notation.
Let $A_i,B_i\subset\E(G_1)$ consist of the edges colored by \one\/
and \two\/ respectively up to the $i$-th round of $\gamesym{G_1}$ and
$A'_i,B'_i\subset\E(G_0)$ consist of the edges colored by \oneo\/
and \twoo\/ respectively up to the $i$-th round of the simulated
game $\gamesym{G_0}$. Initially $A_0=B_0=A'_0=B'_0=\emptyset$.
It will be the case that, up to the $(2i-1)$-th round of the simulated
game $\efgame{\lineg{G_0},\lineg{G_1}}$, the \oneef\/ and the \twoef\/
choose exactly the vertices in $A_i\cup A'_i\cup B_{i-1}\cup B'_{i-1}$
and, up to the $(2i)$-th round, they choose the vertices in
$A_i\cup A'_i\cup B_i\cup B'_i$.

Assume that $S_0$ succeeds in $i$ rounds of $\gamesym{G_0}$
and $D$ succeeds in $2i$ rounds of $\efgame{\lineg{G_0},\lineg{G_1}}$
irrespective of the other players's strategies. Under this assumption,
we describe the move of \two\/ in the $i$-th round of $\gamesym{G_1}$
and then show that this move is successful.

Assume that \one\/ colors an edge $a$ and hence $A_i=A_{i-1}\cup\{a\}$.
Simulating the $(2i-1)$-th round of $\efgame{\lineg{G_0},\lineg{G_1}}$,
the player \two\/ makes the \oneef\/ choose $a$, a vertex in $\lineg{G_1}$,
and then makes the \twoef\/ apply the strategy $D$. Let $a'$ denote the
vertex chosen by the \twoef\/ in $\lineg{G_0}$.
Simulating the $i$-th round of $\gamesym{G_0}$, the player \two\/
makes \oneo\/ color the edge $a'$ thereby setting $A'_i=A'_{i-1}\cup\{a'\}$
and then makes \twoo\/ apply the strategy $S_0$. Let $b'$ denote the
edge colored by \twoo\/ and $B'_i=B'_{i-1}\cup\{b'\}$. Next \two\/
simulates the $(2i)$-th round of $\efgame{\lineg{G_0},\lineg{G_1}}$.
He makes the \oneef\/ choose $b'$, a vertex in $\lineg{G_0}$,
and then makes the \twoef\/ apply $D$. Let $b$ denote the
vertex chosen by the \twoef\/ in $\lineg{G_1}$. Finally, \two\/ colors
the edge $b$ and hence $B_i=B_{i-1}\cup\{b\}$.

We now have to show that $S(D,S_0)$ succeeds in the $i$-th round
irrespective of the player \one's strategy.
Since by our assumption $S_0$ succeeds against any strategy
of \oneo, the subgraphs $A'_i$ and $B'_i$ of $G_0$ are isomorphic.
By the definition of a line graph, the subgraphs of $\lineg{G_0}$
induced by the vertex sets $A'_i$ and $B'_i$ are isomorphic too.
As easily seen from the description of $S(D,S_0)$, the \twoef\/
constructs $A'_i$ from $A_i$ and $B_i$ from $B'_i$.
Since by our assumption $D$ succeeds against any strategy
of the \oneef, the subgraphs induced by $A_i$ in $\lineg{G_1}$ and
by $A'_i$ in $\lineg{G_0}$ are isomorphic, as well as
the subgraphs induced by $B_i$ and $B'_i$ are isomorphic. It follows that
the subgraphs of $\lineg{G_1}$ induced by $A_i$ and $B_i$ are isomorphic.
Since these are the line graphs of the subgraphs $A_i$ and $B_i$
of $G_1$, the latter two are edge-isomorphic. By the condition
imposed on $G_1$, neither $A_i$ nor $B_i$ have a connected component
$K_3$. By the Whitney theorem, we conclude that $A_i$ and $B_i$
are isomorphic and therefore the strategy $S(D,S_0)$ of \two\/
does succeed in $\gamesym{G_1}$ independently of \one's strategy.
\end{proof}

\begin{remark}\label{rem:ext}
Propositions \ref{prop:sym_ef} and \ref{prop:new} actually hold true
for any edge-coloring game in place of $\gamesym G$ if this game has
isomorphism-invariant winning conditions.
\end{remark}

\noindent We are now prepared to prove the claim 1 of Theorem~\ref{thm:odd}.

\begin{corollary}
For odd $n$, $\lsym{P_n}>\frac12\log(n-1)-1$ and
$\lsym{C_n}>\frac12\log n-\frac12$.
\end{corollary}

\begin{proof}
By Proposition \ref{prop:auto}, $\lsym{P_{n+1}}=(n+1)/2$.
A weaker bound $\lsym{P_n}>\frac14\log n - \frac12$
follows from Proposition \ref{prop:sym_ef},
with $G_1=P_n$ and $G_0=P_{n+1}$, and from Proposition \ref{prop:lef}.
To obtain the bound claimed, notice that $\lineg{P_m}=P_{m-1}$
and apply Proposition \ref{prop:new} instead of
Proposition \ref{prop:sym_ef}.
For cycles the proof is the same and uses the fact that $\lineg{C_m}=C_m$.
\end{proof}

\subsection{Upper bound}

We now prove the claim 2 of Theorem \ref{thm:odd}.

\begin{proposition}\label{prop:upper}
If $n$ is odd,
$$
\lsym{P_n}\le(3.5+o(1))\log^2n\mbox{\ and\ }
\lsym{C_n}\le(3.5+o(1))\log^2n.
$$
\end{proposition}

\begin{proof}
We prove the proposition for paths in full detail and then
briefly notice what should be changed for cycles.

Given a subgraph $A$ of $P_n$, we denote its size by $|A|$.
The distance between two subgraphs $A$ and $B$, denoted by
$d(A,B)$, is the minimum distance between vertices $u\in\V(A)$
and $v\in\V(B)$.

We describe a strategy of \one\/ that aims to destroy the
isomorphism between the red and the blue subgraphs possibly sooner.
All moves of \one\/ are split into consecutive
series. The first move of each series creates a new component
of the red subgraph and every subsequent move of the series
prolongs the component in one edge. The component created by \one\/
during the $j$-th series will be denoted by $A_j$.
It will be always the case that
\begin{equation}\label{eq:alpha}
|A_{j+1}|<|A_{j}|.
\end{equation}

\convention{%
Throughout our description of \one's strategy, we assume that
\two\/ plays optimally against this strategy, that is,
keeps the isomorphism between the red and the blue subgraphs
as long as possible. This, together with the condition \refeq{eq:alpha},
implies that in the first move of every series, \two\/ also
must start constructing a new component of the blue subgraph
and in each subsequent move of the series he must extend this
component (or otherwise \two\/ violates the isomorphism and loses
immediately). The component created by \two\/
during the $j$-th series will be denoted by $A'_j$.}

\begin{definition}\label{def:distpair}
In any position of the game such that it is \one's turn, we call
two red components $A_i$ and $A_j$ a {\em distinctive pair\/} if
\begin{enumerate}
\item
$d(A_i,A_j)\ne 2$,
\item
no edge between $A_i$ and $A_j$ has been chosen by the players,
\item
$d(A'_i,A'_j)\ne d(A_i,A_j)$ or between $A'_i$ and $A'_j$
there is at least one edge chosen by the players.
\end{enumerate}
\end{definition}

\notation{%
We number all edges of $P_n$ from one end edge to the other end edge.
For notational convenience we identify the edges with their
numbers $1,2,\ldots,n$.
By $s_j$ and $f_j$ ($s'_j$ and $f'_j$ resp.) we denote the edges
chosen by \one\/ (\two\/ resp.) in the first and the last moves
of the $j$-th series. Note that it is unnecessary that $s_j$ and $f_j$
are the end edges of $A_j$ but most often this will be so.}

Set
\begin{equation}\label{eq:t}
t=4\lceil\log n\rceil+22.
\end{equation}
With the exception of a few last series,
the number of moves in the $j$-th series and hence the length of $A_j$
will be $t-j$. The parameter $t$ is chosen large enough to ensure that,
until the end of the game, $t-j$ is a positive number;
the proof is given by Claim 1 below.

To avoid separately handling several exceptional cases
of small $n$, we just assume $n$ to be sufficiently large
to satisfy the inequality
\begin{equation}\label{eq:n}
n>14t.
\end{equation}
The cases of smaller $n$ are covered by the $o(1)$ term
in the statement of Proposition~\ref{prop:upper}.

\smallskip

\noindent
{\sc The first series of moves by \one.}
In the first move \one\/ chooses the middle edge of $P_n$,
that is, $s_1=(n+1)/2$. Without loss of generality assume that
$s'_1<s_1$. Then in the next moves \one\/ chooses the edges
$s_1+1,s_1+2,\ldots,s_1+t-2=f_1$.

\smallskip

In our further description of \one's strategy, we distinguish
two phases of the game.

\phase 1: \one\/ enforces appearance of a distinctive pair $\{C_0,D_0\}$.

\smallskip

\noindent
{\sc The $j$-th series of moves, $j>1$.}
We assume that after completion of the preceding series the
following conditions are met for $m=j-1$.

\condition 0
$|A_1|>|A_2|>\cdots|A_{m-1}|>|A_m|$.

\condition 1
$n-f_m>2$.

\condition 2
No vertex on the right from $f_m$ has been chosen by the players.

\condition 3
The edges of $A_1,A_2,\ldots,A_m$ have been chosen in the ascending
order. In particular,
$$
s_1<f_1<s_2<f_2<\cdots<s_m<f_m
$$
and $s_p$ and $f_p$ are the end edges of $A_p$, for every $p\le m$.

\condition 4
For notational convenience, define $A_0$ to be the subgraph
of $P_n$ induced by the first vertex of $P_n$ ($A_0$ has a single vertex
and no edge). Let $q$ be such that $A'_m$ is between $A_{q-1}$ and $A_q$.
Then $A'_q,A'_{q+1},\ldots,A'_{m-1}$ all are also between
$A_{q-1}$ and $A_q$, exactly in this order (in the direction
either from $A_{q-1}$ to $A_q$ or from $A_{q}$ to $A_{q-1}$).
In addition, $d(A'_{p-1},A'_p)=d(A_{p-1},A_p)$ for every $q<p\le m$.

\condition 5
$n-(s_p-1)>d(A_{p-1},A_p)$ for every $1<p\le m$.

\smallskip

Observe that Conditions 0--5 are obeyed after the first series of moves,
that is, they are true for $m=1$. Condition 1 follows from \refeq{eq:n};
Conditions 2 and 3 follow from the description of the first series;
Conditions 0, 4, and 5 for the first series are trivial. The fulfillment
of Conditions 0--5 for every series excepting the last series of Phase 1
will be proven in Claim 1 below.

Define a function $\phi$ by $\phi(x)=x+\lceil(n-x)/2\rceil$.
In the first move of the $j$-th series \one\/ chooses
\begin{equation}\label{eq:s_j}
s_j=\phi(f_{j-1}).
\end{equation}
\comment 1
{
This choice of $s_j$ implies Condition 5 for $m=j$, unless the $j$-th
series concludes Phase 1. The function $\phi$ provides the smallest value
of $s_j$ with this property. Note also that $s_j$ really starts
a new component due to Condition 1 for $m=j-1$.
}

The further moves of \one\/ depend much on the first move of \two\/
in the series.

\case 1 {$s'_j<s_j$}%
\one\/ continues the series choosing $s_j+1,s_j+2,\ldots,f_j$.
The last edge in the series is determined from the following rules.
\rrule 1
If $n-(s_j-1)\le t-j$, then $f_j=n$. Otherwise:
\rrule 2
If $n-(\phi(s_j+t-j-1)-1)\ge t-(j+1)$, then $f_j=s_j+t-j-1$.
\rrule 3
If $n-(\phi(s_j+t-j-1)-1)< t-(j+1)$, then
$f_j$ is the smallest number such that $n-\phi(f_j)<f_j-s_j$.
\comment 2
{
These rules can be reformulated as follows. $A_j$ is constructed
edge by edge in the ascending order starting from $s_j$ so that
$|A_j|=t-j$ with two exceptional cases:
{
\renewcommand{\theenumi}{\roman{enumi}}
\renewcommand{\labelenumi}{(\theenumi)}
\begin{enumerate}
\item
Assignment $|A_j|=t-j$ and starting the next component $A_{j+1}$
from $s_{j+1}=\phi(f_j)$ in the ascending order could not give
$|A_{j+1}|=t-(j+1)$ because the final edge $n$ would be reached
earlier than in $t-(j+1)$ moves. In this case $|A_j|$ is taken as smaller
than $t-j$ as possible to keep the relation $|A_{j+1}|<|A_j|$,
where $A_{j+1}$ starts at $\phi(f_j)$ and finishes at $n$.
\item
$|A_j|$ is shorter because the last edge $n$ is reached.
\end{enumerate}
}

\noindent
The case (i) corresponds to Rule 3, and the case (ii) corresponds
to Rule 1.

It is also useful to make the following observation.
Assume that Case 1 occurs in the $(j-1)$-th and in the $j$-th series.
Then, if the case (i) occurs in the $(j-1)$-th series,
the case (ii) must occur in the $j$-th series. Vice versa, if
the case (ii) occurs in the $j$-th series, then the case (i)
must occur in the $(j-1)$-th series.
}

\subcase {1-a} {$f_{j-1}<s'_j<s_j$}%
If $f_j<n$, then Phase 1 of the game continues and \one\/ starts
the $(j+1)$-th series. If $f_j=n$, then \two\/ by \refeq{eq:s_j}
has not enough room to make $A'_j$ so long as $A_j$. Therefore
the isomorphism is violated and \two\/ loses.

\subcase {1-b} {$s'_j<s_{j-1}$}%
Let $q$ be as in Condition 4 for $m=j-1$.
If $A'_j$ is not between $A_{q-1}$ and $A_q$, then
$A_{j-1}$ and $A_j$ are a distinctive pair. The items 2 and 3
of Definition \ref{def:distpair} are clear, and the item 1
is proved in Claim 2 below.
\one\/ therefore terminates Phase 1 and takes $C_0=A_{j-1}$, $D_0=A_j$.
Suppose that $A'_j$ is between $A_{q-1}$ and $A_q$.

If $d(A'_j,A'_{j-1})\ne d(A_j,A_{j-1})$, then
$A_{j-1}$ and $A_j$ again are a distinctive pair and \one\/
terminates Phase 1 with $C_0=A_{j-1}$, $D_0=A_j$.
Suppose that $d(A'_j,A'_{j-1})=d(A_j,A_{j-1})$.
If Condition 4 is violated for $m=j$, then $A'_j$ must be between
$A'_{j-1}$ and $A'_{j-2}$ and therefore  $A_{j-1}$ and $A_{j-2}$
become a distinctive pair. In this case \one\/
terminates Phase 1 and takes $C_0=A_{j-1}$, $D_0=A_{j-2}$.
A distinctive pair does not exist in the only case that
Condition 4 holds true for $m=j$. In this case Phase 1 of the
game continues and \one\/ starts the $(j+1)$-th series of moves.

Let us pay special attention to
the case when $f_j=n$. Condition 4 then cannot happen
in view of Condition 5 for $p=q$. Therefore, a distinctive pair
exists and the game goes to Phase 2.

\case 2 {$s'_j>s_j$}%
\one\/ continues the series choosing
$s_j-1,s_j-2,\ldots,f_j=\max\{f_{j-1}+2,s_j-t+j+1\}$ unless
this maximum equals $f_{j-1}+3$. In the latter case the series is shorter
in one move and $f_j=f_{j-1}+4$. The series
$A_{j-1}$ and $A_j$ are a distinctive pair and
\one\/ terminates Phase 1 with $C_0=A_{j-1}$, $D_0=A_j$.

\comment 3
{
In other words, $A_j$ is constructed starting from $s_j$,
edge by edge in the descending order, until $|A_j|=t-j$
or $d(A_j,A_{j-1})=1$. Special care is taken to ensure that
$d(A_j,A_{j-1})\ne 2$, one of the defining properties of a distinctive pair.
}

\nopagebreak

\noindent
{\sc End of description of Phase 1}

\bigskip

The above program for \one\/ makes sense as long as $t-j$, the length
assigned to $A_j$, is a positive number, which is the case
for at most $t-1$ series of moves. In fact,
we prove that $\lceil\log n\rceil$ series suffice
for \one\/ to terminate Phase 1, that is, either to win the game
or to find a pair $\{C_0,D_0\}$ (Item 3 of Claim 1 below).
We will prove that the pair $\{C_0,D_0\}$ found by \one\/ is indeed
distinctive (Claim 2). We also should prove our assumption that
Conditions 0--5 hold true at the start of every series of moves
within Phase 1
(Claim 1, Item 6). To verify \refeq{eq:alpha}, in addition to
Condition 0 we need to prove the inequality $|A_j|<|A_{j-1}|$
for components created in the last two series of Phase 1 (Claim 1, Item 4).

\begin{claim}{1}
Let \one\/ play as described above and \two\/ play optimally against
this strategy of \one. Suppose that \one\/ has made the $j$-th series of
moves in Phase 1. Then
\begin{enumerate}
\item
$d(A_j,A_{j-1})<\frac12 d(A_{j-1},A_{j-2})$, if $j\ge 3$.
\item
$d(A_j,A_{j-1})<(\frac12)^{j-1} n$, if $j\ge 2$.
\item
$j\le\lceil\log n\rceil$.
\item
$|A_j|<|A_{j-1}|$, if $j\ge 2$.
\item
$|A_j|>\log n+4$.
\item
If the $j$-th series is not last in Phase 1, then Conditions 0--5
hold true for $m=j$.
\end{enumerate}
\end{claim}

\begin{subproof}
We proceed by induction on $j$. Consider two base cases $j=1,2$.
Item 1 is trivial. Item 2 reads $d(A_2,A_1)<\frac12n$ and is
straightforward by the description of \one's strategy. Item 3
is equivalent to $n>2$ and follows from \refeq{eq:n} and \refeq{eq:t}.
By \refeq{eq:n} and by the description of \one's strategy,
$|A_1|=t-1$ and $|A_2|=t-2$ and hence Items 4 and 5 are true.
Taking into account \refeq{eq:n}, it is also easy to check Item 6.

Assume that Items 1--6 are true for the $(j-1)$-th series and
prove each of them for the $j$-th series.

\itemm 1
By the induction assumption applied to Item 5 we have
\begin{equation}\label{eq:b8}
|A_{j-1}|>4.
\end{equation}

Assume that $A_j$ is created in Case 1. Then
\begin{eqnarray*}
&s_{j-2}<f_{j-2}<s_{j-1}<f_{j-1}<s_j<f_j,&\\
&d(A_j,A_{j-1})=s_j-f_{j-1}-1,
\quad
d(A_{j-1},A_{j-2})=s_{j-1}-f_{j-2}-1.&
\end{eqnarray*}
By the choice of $s_{j-1}$ (see \refeq{eq:s_j}),
$$
d(A_{j-2},A_{j-1})+2\ge |A_{j-1}|+d(A_j,A_{j-1})+|A_j|+(n-f_j).
$$
By the choice of $s_{j}$,
$$
|A_j|+(n-f_j)=n-s_j+1>d(A_j,A_{j-1}).
$$
Taking into account \refeq{eq:b8}, we infer that
\begin{equation}\label{eq:b9}
d(A_{j-2},A_{j-1})>2d(A_j,A_{j-1}).
\end{equation}

Assume now that $A_j$ is created in Case 2. Then
\begin{eqnarray*}
&s_{j-2}<f_{j-2}<s_{j-1}<f_{j-1}<f_j<s_j,&\\
&d(A_j,A_{j-1})=f_j-f_{j-1}-1,
\quad
d(A_{j-1},A_{j-2})=s_{j-1}-f_{j-2}-1.&
\end{eqnarray*}
By the choice of $s_{j-1}$,
$$
d(A_{j-2},A_{j-1})+2\ge |A_{j-1}|+d(A_j,A_{j-1})+|A_j|+(n-s_j).
$$
By the choice of $s_{j}$,
$$
n-s_j+1>d(A_j,A_{j-1})+(|A_j|-1).
$$
Taking into account \refeq{eq:b8}, we again easily infer~\refeq{eq:b9}.

\itemm 2
By the induction assumption,
$$
d(A_{j-1},A_{j-2})<(\frac12)^{j-2}n.
$$
Using Item 1, we derive
$$
d(A_j,A_{j-1})<\frac12d(A_{j-1},A_{j-2})<(\frac12)^{j-1}n
$$
as required.

\itemm 3
As $d(A_j,A_{j-1})\ge 1$, this is a consequence of Item 2.

\itemm 4
Note that the component $A_{j-1}$ followed by $A_j$ can be
constructed according to one of six scenarios.
\begin{center}
\begin{tabular}{llll}
Scenario 1: & Case 1, Rule 2 & followed by & Case 1, Rule 2. \\
Scenario 2: & Case 1, Rule 2 & followed by & Case 1, Rule 1. \\
Scenario 3: & Case 1, Rule 2 & followed by & Case 1, Rule 3. \\
Scenario 4: & Case 1, Rule 2 & followed by & Case 2. \\
Scenario 5: & Case 1, Rule 3 & followed by & Case 1, Rule 1. \\
Scenario 6: & Case 1, Rule 3 & followed by & Case 2.
\end{tabular}
\end{center}
In Scenarios 1--4 we have $|A_{j-1}|=t-(j-1)$ and $|A_j|\le t-j$.
In Scenario 5, $|A_j|=n-(s_j-1)$ by Rule 2 and the inequality
$|A_j|<|A_{j-1}|$ is enforced by Rule 3.
In Scenario 6, $|A_j|$ is even shorter than in Scenario 5 because
$|A_j|\le s_j-f_{j-1}-1$
and $s_j-f_{j-1}-1<n-(s_j-1)$ by the choice of $s_j$.

\itemm 5
We distinguish the same six scenarios as above.

\smallskip

\noindent
Scenarios 1--2. We have $|A_j|=t-j$ and the claim follows from Item 3
and \refeq{eq:t}.

\smallskip

\noindent
Scenario 3 will be considered a bit later.

\smallskip

\noindent
Scenario 4. We have $|A_{j-1}|=t-(j-1)$. Since $A_{j-1}$ is constructed
according to Rule 2, we have $n-(s_j-1)\ge t-j$.
Together with \refeq{eq:s_j}, this implies that
$s_j-f_{j-1}-1\ge t-j-2$. As in Case 2
$|A_j|\ge\min\{t-j-1,s_j-f_{j-1}-1\}$, we have $|A_j|\ge t-j-2$.
The claim now follows from Item 3 and \refeq{eq:t}.

\smallskip

\noindent
Scenario 5.
According to Rule 3,
\begin{equation}\label{eq:b1}
n-\phi(f_{j-1})<f_{j-1}-s_{j-1}
\end{equation}
and
\begin{equation}\label{eq:b2}
n-\phi(f_{j-1}-1)\ge (f_{j-1}-1)-s_{j-1}.
\end{equation}
{}From \refeq{eq:b1} we infer
\begin{equation}\label{eq:b3}
n-s_{j-1}<3(f_{j-1}-s_{j-1})+1,
\end{equation}
and from \refeq{eq:b2} we infer
\begin{equation}\label{eq:b4}
n-f_{j-1}+1\ge 2(f_{j-1}-s_{j-1}-1).
\end{equation}
As $A_{j-1}$ is constructed according to Rule 3 and therefore
the assumption of Rule 1 is false, we have
\begin{equation}\label{eq:b5}
n-(s_{j-1}-1)>t-(j-1).
\end{equation}
{}From \refeq{eq:b3} and \refeq{eq:b5} we conclude that
\begin{equation}\label{eq:b6}
f_{j-1}-s_{j-1}>(t-j-1)/3.
\end{equation}
By Rule 1, $|A_j|$ is equal to
$$
n-(s_j-1)=n-\phi(f_{j-1})+1\ge(n-f_{j-1}+1)/2
$$
Using \refeq{eq:b4} and \refeq{eq:b6}, we obtain
\begin{equation}\label{eq:b7}
n-(s_j-1)\ge f_{j-1}-s_{j-1}-1>(t-j-1)/3-1.
\end{equation}
Hence $|A_j|>(t-j-4)/3$ and the claim follows from Item 3 and \refeq{eq:t}.

\smallskip

\noindent
Scenario 6.
Since $A_{j-1}$ is constructed according to Rule 3, we have
$$
(s_j-1)-f_{j-1}\le n-(s_j-1)-1\le (t-j)-2.
$$
In Case 2 we have $|A_j|\ge\min\{t-j-1,s_j-f_{j-1}-1\}$ and hence
$|A_j|\ge s_j-f_{j-1}-1$. The latter value, by the choice of $s_j$,
is no less than $n-(s_j-1)-2$.
Similarly to Scenario 5, the relation \refeq{eq:b7} is true and hence
$$
|A_j|\ge n-(s_j-1)-2>(t-j-10)/3.
$$
It remains to apply Item 3 and \refeq{eq:t}.

\smallskip

\noindent
Scenario 3.
By Rule 3, $|A_j|>n-(s_{j+1}-1)$. Applying precisely the same
argument as in Scenario 5, similarly to \refeq{eq:b7} we derive
$$
n-(s_{j+1}-1) > (t-(j+1)-4)/3.
$$
It remains to apply Item 3 and \refeq{eq:t}.

\itemm 6
The assumption made in this item implies that for the $j$-th series
we have Case 1 and $f_j<n$, that is, either Rule 2 or Rule 3 was applied.
Condition 0 follows from the induction assumption and Item 4.
Condition 1 is ensured by Rules 2 and 3.
Conditions 2 and 3 can be violated only in Case 2 which always terminates
Phase 1.
Condition 4 is obvious if $q=m$. If $q\le m-1$,
the condition follows from the induction assumption
(see explanations accompanying the description of Case 1).
Condition 5 follows from the choice of $s_p$ (see \refeq{eq:s_j})
and Condition 3.
\end{subproof}

\begin{claim}{2}
If \one\/ finishes Phase 1 with some $C_0$ and $D_0$, then these
components are a distinctive pair.
\end{claim}

\begin{subproof}
If $l$ is the number of series in Phase 1, then either
$\{C_0,D_0\}=\{A_{l},A_{l-1}\}$ or $\{C_0,D_0\}=\{A_{l-1},A_{l-2}\}$.
It is easy to check that this pair is always chosen so that
the items 2 and 3 of Definition \ref{def:distpair} are true.
Let us check the item 1, i.e., $d(C_0,D_0)\ne 2$.

If $\{C_0,D_0\}$ is chosen in Case 1, then $d(C_0,D_0)$
equals either $s_l-f_{l-1}-1$ or $s_{l-1}-f_{l-2}-1$.
By the choice \refeq{eq:s_j} of $s_j$ these values are not less
than $|A_l|-2$ and $|A_{l-1}|-2$ respectively. By Item 5 of Claim 1,
$d(C_0,D_0)>2$.

If $\{C_0,D_0\}=\{A_l,A_{l-1}\}$ is chosen in Case 2, then the inequality
$d(C_0,D_0)\allowbreak\ne 2$ is true by the choice of the length $|A_l|$ in Case 2.
\end{subproof}

\notation{%
In the sequel we denote the number of series in Phase 1 by $l$.
Let $t'=|A_l|$, the length of the shortest component created in Phase 1.}

\begin{claim}{3}
\mbox{}
\begin{enumerate}
\item
$l\le\lceil\log n\rceil$.
\item
$t'>\log n + 4$.
\end{enumerate}
\end{claim}

\begin{subproof}
The claim is a direct corollary of Items 3 and 5 of Claim 1.
\end{subproof}

\phase 2: \one\/ reduces the distance between components of a
distinctive pair to~1.

\smallskip

\noindent
{\sc The $(l+j)$-th series of moves (the $j$-th series in Phase 2).}
Let $\{C_{j-1},D_{j-1}\}$ be the distinctive pair created in the preceding
series. In particular, $\{C_0,D_0\}$ is the output of Phase 1.
If $d(C_{j-1},D_{j-1})=1$, \one\/ makes the last move described below.
Suppose that $s(C_{j-1},D_{j-1})\ge 3$. Let $a$ and $b$ be two nearest
edges in $C_{j-1}$ and $D_{j-1}$ respectively.
Without loss of generality, assume
$a<b$. In the first move of the series \one\/ chooses the medium edge
$s_{l+j}=\lceil(a+b-1)/2\rceil$. The further moves of \one\/ depend on
the first move of~\two.

\case 1 {$a<s'_{l+j}<s_{l+j}$}%
\one\/ continues the series choosing
$s_{l+j}+1,s_{l+j}+2,\ldots,f_{l+j}=\min\{b-2,s_{l+j}+(t'-j-1)\}$
unless $s_{l+j}+(t'-j-1)=b-3$. In the latter case \one\/ stops
one move earlier at $f_{l+j}=b-4$. The new distinctive pair is
$C_j=A_{l+j}$ and $D_j=D_{j-1}$.

\case 2 {$s_{l+j}<s'_{l+j}<b$}%
\one\/ continues the series choosing
$s_{l+j}-1,s_{l+j}-2,\ldots,f_{l+j}=\max\{a+2,s_{l+j}-(t'-j-1)\}$
unless $s_{l+j}-(t'-j-1)=a+3$. In the latter case \one\/ stops
one move earlier at $f_{l+j}=a+4$. The new distinctive pair is
$C_j=C_{j-1}$ and $D_j=A_{l+j}$.

\case 3 {$s'_{l+j}$ is not between $C_{j-1}$ and $D_{j-1}$}%
\one\/ continues the series choosing
$s_{l+j}+1,s_{l+j}-1,s_{l+j}+2,s_{l+j}-2$ and so on until one
of the following situations happens:
\begin{enumerate}
\item
$d(C_{j-1},A_{l+j})=d(A_{l+j},D_{j-1})=1$.
\item
$|A_{l+j}|=t'-j$ but $d(C_{j-1},A_{l+j})\ne2$ and $d(A_{l+j},D_{j-1})\ne2$.
\item
$t'-j-3\le|A_{l+j}| < t'-j$ and
$d(C_{j-1},A_{l+j})=d(A_{l+j},D_{j-1})=3$.
\end{enumerate}
As shown in Claim 5 below, at least one of the pairs
$\{C_{j-1},A_{l+j}\}$ or $\{A_{l+j},D_{j-1}\}$ is distinctive and
\one\/ takes it as $\{C_j,D_j\}$.

\smallskip

\noindent
{\sc The last move.}
As soon as \one\/ creates a distinctive pair $\{C,D\}$ with $d(C,D)=1$,
he chooses the edge between $C$ and $D$ and wins. To keep isomorphism,
\two\/ should make, in place of two corresponding blue components
$C'$ and $D'$, a new component of length $|C'|+|D'|+1$. This task
is impossible to implement as $d(C',D')>1$.

\medskip

The description of \one's strategy in Phase 2 makes sense as long as
the value assigned to the length of a series is a positive number.
The smallest value that can be assigned for the $(l+j)$-th series,
if it is not last in Phase 2, is $t'-j-3$. Hence the
condition $j<t'-4$ is required. The same condition is required
also in order to pass by the forbidden distance $d(C_j,D_j)=2$
in cases $d(C_{j-1},D_{j-1})=4,5,6$.
We will show that $t'-4$ series are indeed enough for \one\/
to finish Phase 2 (Claim 4, Item 3) and that $\{C,D\}$, the outcome
of Phase 2, is indeed a distinctive pair with $d(C,D)=1$ (Claim 5).
We also will verify \refeq{eq:alpha} for the components created
during Phase 2 (Claim 4, Item 4).

\begin{claim}{4}
Let \one\/ play as described above and let \two\/ play optimally
against this strategy of \one. Suppose that \one\/ has made
the $(l+j)$-th series of moves in Phase 2. Then
\begin{enumerate}
\item
$d(C_j,D_j)\le\frac12d(C_{j-1},D_{j-1})$.
\item
$d(C_j,D_j)\le(\frac12)^jd(C_0,D_0)$.
\item
$j<\log n-1$.
\item
$|A_{l+j}|<|A_{l+(j-1)}|$.
\end{enumerate}
\end{claim}

\begin{subproof}
Item 1 is clear from \one's strategy. Item 2 follows from Item 1.
Item 3 follows from Item 2, because $d(C_j,D_j)\ge1$ and
$d(C_0,D_0)<\frac12n$.

Item 4 is given by easy inspection of \one's
strategy. If the $(l+j)$-th series is neither last nor last but one
in Phase 2, then $|A_{l+(j-1)}|=t'-(j-1)$ and $|A_{l+j}|=t'-j$.
If the $(l+j)$-th series is last but one, then
$|A_{l+(j-1)}|=t'-(j-1)$ and $t'-j-3\le|A_{l+j}|\le t'-j$.
If the $(l+j)$-th series is last, then either
$|A_{l+(j-1)}|=t'-(j-1)$ and $|A_{l+j}|\le t'-j$ or
$t'-(j-1)-3\le|A_{l+(j-1)}|\le t'-(j-1)-1$ and $|A_{l+j}|=1$.
In the latter case Item 4 follows from Item 3 and Claim 3.
\end{subproof}

\begin{claim}{5}
Let \one\/ play as described above and let \two\/ play optimally
against this strategy of \one. Suppose that \one\/ has made
the $(l+j)$-th series of moves in Phase 2. Then
$\{C_j,D_j\}$ is a distinctive pair and eventually $d(C_j,D_j)=1$.
\end{claim}

\begin{subproof}
Conditions 1 and 2 in Definition \ref{def:distpair} are directly enforced
by \one's strategy. Condition 3 is easy to see if the $(l+j)$-th
series of moves was done in Cases 1 or 2. Assume it was done
in Case 3. We have to prove that at least one of the pairs
$\{C_{j-1},A_{l+j}\}$ or $\{A_{l+j},D_{j-1}\}$ is distinctive.
Suppose, to the contrary, that they both are not. This implies for
the corresponding blue components $C',D',A'_{l+j}$ that no edge
is chosen between $C',A'_{l+j}$ and between $A'_{l+j},D'$, and that
$d(C',A'_{l+j})=d(C_{j-1},A_{l+j})$ and $d(A'_{l+j},D')=d(A_{l+j},D_{j-1})$.
It follows that before the $(l+j)$-th series no edge was chosen between $C'$
and $D'$ and $d(C',D')=d(C_{j-1},D_{j-1})$. Thus,
$C_{j-1}$ and $D_{j-1}$ could not be a distinctive pair, a contradiction.
\end{subproof}

It remains to estimate the total number of moves in the game if \one\/
follows the above strategy. Since $\{C_0,D_0\}$ is always either
$\{A_l,A_{l-1}\}$ or $\{A_{l-1},A_{l-2}\}$, we have
$d(C_0,D_0)\le d(A_{l-1},A_{l-2})$. By Item 2 of Claim 4 and
Item 2 of Claim 1, in the last $(l+k)$-th series of moves of Phase 2
we have
$d(C_k,D_k)<(\frac12)^{l+k-2}n$ and therefore $l+k\le\lfloor\log n\rfloor+2$.
As the $j$-th series has at most $t-j$ moves, the total number
of rounds in the game does not exceed
$\sum_{j=1}^{\lfloor\log n\rfloor+2}(t-j)+1=3.5\log^2 n+O(\log n)$.

The proof of Proposition \ref{prop:upper} for paths is complete.
\end{proof}


\medskip

\noindent
{\sl Proof-sketch of Proposition \ref{prop:upper} for cycles.}
We employ the same idea as for paths and refer to strategies in
Phases 1 and 2 described above. The moves of \one\/ are split
into series, and in a series \one\/ creates a component of the red
subgraph. The $j$-th series consists of $t-j$ moves, with a few possible
exceptions in the end of the game.

We adopt the notion of a distinctive pair of components with
the only refinement: The distance $d(A,B)$ between two components
$A$ and $B$ is the minimum length of a path that joins a vertex
in $A$ and a vertex in $B$ and that consists of edges unchosen so far.
Thus, $d(A,B)$ may differ from the standard distance in a graph.
The goal of \one\/ is to create a distinctive pair and then to apply
the strategy of Phase 2 literally. However, creation of a distinctive pair
in cycles is a bit more complicated task. Namely, before applying
the strategy of Phase 1 some additional efforts are needed.

In the first series of moves \one\/ creates the component $A_1$
and \two, not to lose immediately, creates the component $A'_1$
of the same length. Denote two paths connecting $A_1$ and $A'_1$
by $I_1$ and $I_2$, and their lengths by $l_1$ and $l_2$.
Note that one number of $l_1$ and $l_2$ is odd and the other is even.
Without loss of generality assume $l_1>l_2\ge0$.

\case 1 {$l_1$ is odd}%
\one\/ starts the second series choosing $s_2$, the middle
edge of $I_1$. If \two\/ chooses $s'_2$ in $I_1$ between $s_2$
and $A'_1$, then \one\/ completes $A_2$ and $\{A_1,A_2\}$ is a
distinctive pair. If \two\/ chooses $s'_2$ in $I_1$ between
$A_1$ and $s_2$, then \one\/ continues to play on $I_1$ in
the direction towards $A'_1$ applying the strategy of Phase 1.

\two\/ has another possibility to try to avoid creating a
distinctive pair: He can choose $s'_2$ in $I_2$ at the same
distance from $A'_1$ as between $A_1$ and $s_2$. In this case,
if \one\/ continues to play Phase 1 on $I_1$ in
the direction towards $A'_1$, \two\/ can copy
moves of \one\/ in $I_2$. Nevertheless, since $l_2<l_1$, eventually
either the isomorphism will be violated, or a distinctive pair appears,
or \two\/ will be forced to switch back to $I_1$.

\case 2 {$l_1$ is even}%
Playing on $I_1$ gives no gain for \one\/ because \two\/ can keep
isomorphism using the involutory fixed-edge-free automorphism of $I_1$
(this is a difference between the cases of paths and cycles).
Therefore, \one\/ should play in $I_2$. However, it is impossible
for \one\/ to adapt the strategy of Phase 1 directly because
\two\/ can just copy moves of \one\/ in $I_1$. To prevent this,
\one\/ chooses $s_2$ in $I_1$ at distance $(l_2-1)/2$ from $A_1$.

If \two\/ chooses $s'_2$ in $I_1$ at the same distance from $A'_1$,
then \one\/ completes $A_2$ and starts the third series choosing
$s_3$ at the center of $I_2$. After this everything goes through
as in Case 1 with roles of $I_1$ and $I_2$ interchanged. Note
that \two\/ is not able to choose $s'_3$ in $I_1$ at the same
distance from $A'_1$ as between $s_3$ and $A_1$ because the corresponding
edge is already occupied in the preceding series.

If \two\/ chooses $s'_2$ in $I_1$ between $s_2$ and $A'_1$ but
not at distance $(l_2-1)/2$ from $A'_1$, then \one\/ completes
$A_2$ so that $\{A_1,A_2\}$ is a distinctive pair.

If \two\/ chooses $s'_2$ in $I_1$ between $A_1$ and $s_2$
or chooses $s'_2$ to be the middle edge of $I_2$, then
\one\/ continues the game on $I_1$ as in Case 1.

\begin{remark}
Notice an essential difference between the \EF\/ game and the
symmetry breaking-preserving game. While in each round of the
former game \two\/ is obliged to extend the isomorphism established
in the preceding round, in the latter game no dependence between
isomorphisms in two successive rounds is required. Eliminating this
difference, let $\symplus G$ be a modification of the symmetry
breaking-preserving game $\gamesym G$ in which \two\/ not merely
keeps the red and the blue subgraphs isomorphic but, moreover,
extends the isomorphism between them from round to round.
Clearly, $\lsymplus G\le \lsym G$.

Propositions \ref{prop:sym_ef}
and \ref{prop:new} hold true for $\symplus G$ with minor
changes in the proofs. In particular, in the proof of Proposition
\ref{prop:new} we need to apply a stronger form of the Whitney
theorem asserting that, with a few exceptions excluded by prohibiting
a subgraph $K_3$, an isomorphism between $\lineg{H_1}$ and
$\lineg{H_2}$ is induced by an isomorphism between $H_1$ and $H_2$
and, moreover, the latter is unique for graphs of size more than 1.
Since Proposition \ref{prop:auto} holds true for $\symplus G$ as well,
we obtain the same lower bounds $\lsymplus{P_n}=\Omega(\log n)$ and
$\lsymplus{C_n}=\Omega(\log n)$ for odd paths and cycles. The upper
bound of Proposition \ref{prop:upper} can be improved to
$\lsymplus{P_n}=O(\log n)$ and $\lsymplus{C_n}=O(\log n)$ for $n$ odd.
The proof becomes much simpler because now the red and blue components
$A_j$ and $A'_j$ correspond to one another by the rules of
the game rather than by having the same distinctive length.
In particular, \one\/ can now make each series of moves being
of constant length.

It would be interesting to know how much the values of $\lsymplus G$
and $\lsym G$ can differ from each other.
\end{remark}

\section{Games on complete graphs}\label{s:k_n}

In this section we analyze the symmetry breaking-preserving game on the
complete graph of order $n$. Unlike the preceding section where we used
the knowledge of $\lsym{P_n}$ and $\lsym{C_n}$ for even $n$, we now have
to estimate $\lsym{K_n}$ for all $n$. The relation with the \EF\/ game
can still give us some information. Similarly to Proposition
\ref{prop:sym_ef} one can prove that, if
$\lef{G_0,G_1}\ge 4\lsym{G_0}+4$, then $\lsym{G_1}\le\lsym{G_0}$.
Since $\lef{K_n,K_{n+1}}=n$, it follows that either $\lsym{K_n}>n/4-1$
for all $n$ or $\lsym{K_n}=O(1)$. We here prove the latter alternative.

\begin{theorem}\label{thm:k_n}
$\lsym{K_n}\le 6$ for all $n$.
\end{theorem}

\begin{proof}
If $n\le 5$, the assertion is trivial.
We assume that $n\ge 6$ and describe a strategy of \one\/
breaking the isomorphism in at most 7 moves. In the first three
rounds \one\/ creates a 3-star in such a way that \two\/ is not able to
choose any edge connecting leafs of this star without immediately losing.
This can be done so that one
of the five positions in Figure \ref{fig:} occurs.

\noindent
\begin{figure}[t]
\centerline{
\unitlength=1.00mm
\special{em:linewidth 0.4pt}
\linethickness{0.4pt}
\begin{picture}(141.00,86.00)
\put(110.00,85.00){\circle*{2.00}}
\put(120.00,85.00){\circle*{2.00}}
\put(130.00,85.00){\circle*{2.00}}
\put(120.00,65.00){\circle*{2.00}}
\bezier{20}(120.00,65.00)(115.00,75.00)(110.00,85.00)
\bezier{20}(120.00,65.00)(119.00,75.00)(120.00,85.00)
\bezier{20}(120.00,65.00)(125.00,75.00)(130.00,85.00)
\put(95.00,40.00){\circle*{2.00}}
\put(105.00,40.00){\circle*{2.00}}
\put(115.00,40.00){\circle*{2.00}}
\put(105.00,20.00){\circle*{2.00}}
\bezier{20}(105.00,20.00)(100.00,30.00)(95.00,40.00)
\bezier{20}(105.00,20.00)(104.00,30.00)(105.00,40.00)
\bezier{20}(105.00,20.00)(110.00,30.00)(115.00,40.00)
\put(130.00,65.00){\circle*{2.00}}
\put(140.00,85.00){\circle*{2.00}}
\emline{120.00}{65.00}{1}{130.00}{65.00}{2}
\emline{130.00}{65.00}{3}{130.00}{85.00}{4}
\emline{130.00}{65.00}{5}{140.00}{85.00}{6}
\emline{140.00}{85.00}{7}{140.00}{85.00}{8}
\put(55.00,85.00){\circle*{2.00}}
\put(65.00,85.00){\circle*{2.00}}
\put(75.00,85.00){\circle*{2.00}}
\put(65.00,65.00){\circle*{2.00}}
\bezier{20}(65.00,65.00)(60.00,75.00)(55.00,85.00)
\bezier{20}(65.00,65.00)(64.00,75.00)(65.00,85.00)
\bezier{20}(65.00,65.00)(70.00,75.00)(75.00,85.00)
\put(75.00,65.00){\circle*{2.00}}
\put(85.00,85.00){\circle*{2.00}}
\emline{65.00}{85.00}{9}{75.00}{65.00}{10}
\emline{75.00}{65.00}{11}{75.00}{85.00}{12}
\emline{75.00}{65.00}{13}{85.00}{85.00}{14}
\put(10.00,85.00){\circle*{2.00}}
\put(20.00,85.00){\circle*{2.00}}
\put(30.00,85.00){\circle*{2.00}}
\put(20.00,65.00){\circle*{2.00}}
\bezier{20}(20.00,65.00)(15.00,75.00)(10.00,85.00)
\bezier{20}(20.00,65.00)(19.00,75.00)(20.00,85.00)
\bezier{20}(20.00,65.00)(25.00,75.00)(30.00,85.00)
\put(30.00,65.00){\circle*{2.00}}
\emline{20.00}{85.00}{15}{30.00}{65.00}{16}
\emline{30.00}{65.00}{17}{20.00}{65.00}{18}
\emline{30.00}{65.00}{19}{30.00}{85.00}{20}
\put(125.00,57.00){\makebox(0,0)[cc]{Position 3}}
\put(70.00,57.00){\makebox(0,0)[cc]{Position 2}}
\put(22.00,57.00){\makebox(0,0)[cc]{Position 1}}
\put(15.00,40.00){\circle*{2.00}}
\put(25.00,40.00){\circle*{2.00}}
\put(35.00,40.00){\circle*{2.00}}
\put(25.00,20.00){\circle*{2.00}}
\bezier{20}(25.00,20.00)(20.00,30.00)(15.00,40.00)
\bezier{20}(25.00,20.00)(24.00,30.00)(25.00,40.00)
\bezier{20}(25.00,20.00)(30.00,30.00)(35.00,40.00)
\put(45.00,40.00){\circle*{2.00}}
\put(55.00,40.00){\circle*{2.00}}
\put(65.00,40.00){\circle*{2.00}}
\emline{45.00}{40.00}{21}{25.00}{20.00}{22}
\emline{25.00}{20.00}{23}{55.00}{40.00}{24}
\emline{25.00}{20.00}{25}{65.00}{40.00}{26}
\put(30.00,8.00){\makebox(0,0)[cc]{Position 4}}
\put(125.00,40.00){\circle*{2.00}}
\put(135.00,40.00){\circle*{2.00}}
\put(125.00,20.00){\circle*{2.00}}
\emline{115.00}{40.00}{27}{125.00}{20.00}{28}
\emline{125.00}{20.00}{29}{125.00}{40.00}{30}
\emline{125.00}{20.00}{31}{135.00}{40.00}{32}
\put(115.00,8.00){\makebox(0,0)[cc]{Position 5}}
\put(18.00,20.00){\makebox(0,0)[lc]{$u_0$}}
\put(15.00,45.00){\makebox(0,0)[ct]{$u_1$}}
\put(25.00,45.00){\makebox(0,0)[ct]{$u_2$}}
\put(35.00,45.00){\makebox(0,0)[ct]{$u_3$}}
\put(45.00,45.00){\makebox(0,0)[ct]{$v_1$}}
\put(55.00,45.00){\makebox(0,0)[ct]{$v_2$}}
\put(65.00,45.00){\makebox(0,0)[ct]{$v_3$}}
\put(95.00,45.00){\makebox(0,0)[ct]{$u_1$}}
\put(105.00,45.00){\makebox(0,0)[ct]{$u_2$}}
\put(115.00,45.00){\makebox(0,0)[ct]{$u_3$}}
\put(125.00,45.00){\makebox(0,0)[ct]{$v_2$}}
\put(135.00,45.00){\makebox(0,0)[ct]{$v_1$}}
\put(105.00,15.00){\makebox(0,0)[cb]{$u_0$}}
\put(125.00,15.00){\makebox(0,0)[cb]{$v_0$}}
\end{picture}
}
\caption{The first three rounds of $\gamesym{K_n}$.
\one's edges are dotted and \two's edges are continuous.}
\label{fig:}
\end{figure}

The next move of \one\/ from Position 1 creates a triangle and
simultaneously blocks creating a triangle by \two.
In Positions 2 and 3 the player \one\/ is able in the next three
moves to create a $K_4$ in such a way that \two\/ cannot do the same.

{\sl Game from Position 4.}
In the next two rounds \one\/ chooses the edges $\{v_1,v_2\}$ and
$\{v_2,v_3\}$. If \two\/ in these rounds chooses two edges of the
triangle $T=\{u_1,u_2,u_3\}$ with the common vertex $u_i$, then
\one\/ chooses $\{u_i,v_2\}$ and wins. Otherwise in the 6-th and 7-th
moves \one\/ chooses two edges of $T$ and wins.

{\sl Game from Position 5.}
In the 4-th and 5-th rounds \one\/ chooses the edges $\{u_3,v_2\}$
and $\{v_2,v_1\}$ respectively. If \two\/ in these rounds chooses
edges not both in $T$, in the next two moves \one\/ chooses two
edges of $T$ and wins. Assume therefore that in the 4-th and 5-th
rounds \two\/ chooses $\{u_3,u_2\}$ and $\{u_2,u_1\}$
(the choice of $\{u_3,u_1\}$ and $\{u_1,u_2\}$ is symmetric).
In the next round \one\/ chooses $\{u_1,v_1\}$ and \two\/ is forced
to choose $\{u_1,v_2\}$. Finally, \one\/ chooses $\{u_2,v_1\}$ and wins.
\end{proof}

\begin{remark}
A more lengthy and complicated analysis allows us to lower the bound 6
of Theorem~\ref{thm:k_n} to 5.
\end{remark}

Finally we briefly discuss the case of
complete bipartite graphs. If at least one of
$m$ and $l$ is even, then $K_{m,l}$ has an involutory automorphism
without fixed edges and, by Proposition \ref{prop:auto},
$\lsym{K_{m,l}}$ is maximum possible for graphs of this size.
If both $m$ and $l$ are odd, $K_{m,l}$ has no involutory
fixed-edge-free automorphism but removal of one edge from $K_{m,l}$
leads to a graph $K_{m,l}-e$ with such an automorphism.
It is therefore interesting to estimate $\lsym{K_{m,l}}$ for
$ml$ odd.

An easy lower bound is
\begin{equation}\label{eq:bipart}
\lsym{K_{m,l}}\ge\max\{\textstyle\frac{m-1}2,\frac{l-1}2\}.
\end{equation}
The appropriate strategy of \two\/ is based on a partial involutory
automorphism of $K_{m,l}$ constructed during the course of the game.
The automorphism leaves one vertex class fixed. Whenever during
the course of the game in the other vertex class a new vertex of
the red subgraph appears, the automorphism interchanges it
with an arbitrary vertex in this class that is unchosen so far.

Note that, if $ml$ is odd, then $K_{m,l}-e$,
$K_{m-1,l}$, $K_{m,l-1}$, and $K_{m-1,l-1}$ all have involutory
automorphisms without fixed edges. One could therefore try to apply
Propositions \ref{prop:sym_ef} and \ref{prop:new} with
$G_1=K_{m,l}$ and $G_0$ one of these graphs. However, in all the
cases $\lef{G_0,G_1}$ and $\lef{\lineg{G_0},\lineg{G_1}}$ are not
large enough to give us anything better than~\refeq{eq:bipart}.

\begin{question}\label{que:}
What are the asymptotics of $\lsym{K_{n,n}}$ for odd $n$?
\end{question}

\bigskip

\noindent
{\sl Note added in proof.}
Question \ref{que:} was recently answered by Oleg Pikhurko who
proved that $\lsym{K_{n,n}}\le 2n+38$ for odd $n\ge 51$. This matches
up to a constant factor the lower bound \refeq{eq:bipart} that reads
$\lsym{K_{n,n}}\ge (n-1)/2$. Pikhurko's result is actually more
general and implies that, if $m\le l\le m^{O(1)}$, then
$\lsym{K_{m,l}}=O(l)$.

\subsection*{Acknowledgement}
We thank Helmut Veith for helpful discussions on the \EF\/ games,
Oleg Pikhurko for useful comments on the paper, and an anonymous referee
for a considerable simplification of our original proof of
Theorem~\ref{thm:k_n}.

\end{document}